\documentclass[fleqn]{mat01}
\usepackage{times,mathtimy,amssymb,latexsym}
\begin{document}

\setcounter{page}{209}
\firstpage{209}

\newtheorem{theore}{Theorem}
\renewcommand\thetheore{\arabic{section}.\arabic{theore}}
\newtheorem{theor}[theore]{\bf Theorem}
\newtheorem{definit}[theore]{\rm DEFINITION}
\newtheorem{lem}[theore]{Lemma}
\newtheorem{rem}[theore]{Remark}
\newtheorem{propo}[theore]{\rm PROPOSITION}
\newtheorem{coro}[theore]{\rm COROLLARY}

\def\remar{\trivlist \item[\hskip \labelsep{\it Remarks.}]}

\title{Some necessary and sufficient conditions for Hypercyclicity
Criterion\thanks{This paper is a part of the second author's Doctoral
thesis, written at Shiraz University under the direction of the first
author.}}

\markboth{B~Yousefi and H~Rezaei}{Hypercyclicity Criterion}

\author{B~YOUSEFI and H~REZAEI}

\address{Department of Mathematics, College of Sciences, Shiraz
University, Shiraz 71454, Iran\\
\noindent E-mail: byousefi@hafez.shirazu.ac.ir}

\volume{115}

\mon{May}

\parts{2}

\pubyear{2005}

\Date{MS received 14 March 2004; revised 29 January 2005}

\begin{abstract}
We give necessary and sufficient conditions for an operator on a
separable Hilbert space to satisfy the hypercyclicity criterion.
\end{abstract}


\maketitle

\section{Introduction}

Suppose that $X$ is a separable topological vector space and $T$
is a continuous linear mapping on $X$. If $x \in X$, then the
orbit of $x$ under $T$ is defined as $Orb(T,x) = \{x, Tx, T^{2}x,
\ldots\}$. An operator $T$ is called hypercyclic if there is a
vector $x$ such that $Orb(T,x)$ is dense in $X$ and in this case
$x$ is called a hypercyclic vector for $T$ (see \cite{14} for an
exhaustive survey on hypercyclicity).

It is interesting that many continuous linear mappings can
actually be hypercyclic. The first example of hypercyclicity
appeared in the space of entire functions, by Birkhoff \cite{3} in 1929. He
showed the hypercyclicity of the translation operator, while
MacLane \cite{19} proved the hypercyclicity of the differentiation
operator in 1952. Hypercyclicity on Banach spaces was discussed in
1969 by Rolewics \cite{20}, who showed that $\lambda B$ is
hypercyclic whenever $B$ is the unilateral backward shift (on
$\ell^{p}$ and $c_{0}$) and $|\lambda| > 1$.

A nice condition for hypercyclicity is the Hypercyclicity
Criterion (Theorem~1.1 below), which was developed by Kitai
\cite{17} and independently by Gethner and Shapiro \cite{12}. This
criterion has been used to show that certain classes of
composition operators \cite{6}, weighted shifts \cite{21},
adjoints of multiplication operators \cite{7}, and adjoints of
subnormal and hyponormal operators \cite{5}, are hypercyclic.
Hypercyclicity has also been established in various other settings
by means of this criterion \cite{1,4,6,8,12,13,16}. Salas
\cite{21} showed that every perturbation of the identity by a
unilateral weighted backward shift with nonzero bounded weights is
hypercyclic, and he also gave a characterization of the
hypercyclic weighted shifts in terms of their weights.
But, then Montes and Leon showed that these hypercyclic operators
do satisfy the criterion as well (\S2 of \cite{17} and
Proposition~4.3 of \cite{18}). Bes and Peris proved that a
continuous linear operator $T$ on a Frechet space satisfies the
Hypercyclicity Criterion if and only if it is hereditarily
hypercyclic. In particular they show that hypercyclic operators
with either a dense generalized kernel or a dense set of periodic
points must satisfy the criterion. Also, they provide a
characterization of those weighted shifts $T$ that are
hereditarily hypercyclic with respect to a given sequence
$\{n_{k}\}_k$ of positive integers, as well as conditions under
which $T$ and $\{T^{n_k}\}_k$ share the same set of hypercyclic
vectors \cite{2}.

\begin{theor}[(The Hypercyclicity Criterion)] Suppose $X$ is a
separable Banach space and $T$ is a continuous linear mapping on
$X$. If there exists two dense subsets $Y$ and $Z$ in $X$ and a
sequence $\{n_{k}\}$ such that{\rm :}
\begin{enumerate}
\renewcommand\labelenumi{\rm \arabic{enumi}.}
\leftskip -.2pc
\item $T^{n_{k}} y \rightarrow 0$ for every $y \in Y${\rm ,}

\item there exists functions $S_{n_{k}}\hbox{\rm :}\ Z \rightarrow X$
such that for every $z \in Z, S_{n_{k}} z \rightarrow 0${\rm ,} and
$T^{n_{k}} S_{n_{k}} z \rightarrow z${\rm ,}\vspace{-.5pc}
\end{enumerate}
then $T$ is hypercyclic.
\end{theor}

Note that the sequence $\{n_{k}\}$ in Theorem~1.1 need not be the
entire sequence $\{n_{k}\} = \{k\}$ of positive integers. Salas
\cite{22} and Herrero \cite{15} have shown that there are
hypercyclic operators on Hilbert spaces that do not satisfy the
Hypercyclicity Criterion for the entire sequence $\{k\}$, but so
far no hypercyclic operator has been found that does not satisfy
the Hypercyclicity Criterion in its general form. In this paper
our work was stimulated by the well-known question: Does every
hypercyclic operator satisfy the hypothesis of the Hypercyclicity
Criterion? (see \cite{2}).

We give necessary and sufficient conditions in terms of open
subsets for an operator on a separable Hilbert space to satisfy
the Hypercyclicity Criterion. For this, see Theorem~2.6,
Corollary~2.11 and Proposition~2.12. Also, in the proof of
Theorem~2.6, we pay attention to hypercyclicity on the operator
algebra $B(H)$ and the algebra of Hilbert--Schmidt operators,
$B_{2}(H)$. Recall that if $\{e_{i}\}_{i}$ is an orthonormal basis
for a separable Hilbert space $H,A \in B(H)$ and
\begin{equation*}
\|A\|_{2} = \left[ \sum^{\infty}_{i=1} \|Ae_{i}\|^{2}
\right]^{1/2},
\end{equation*}
then $\|A\|_{2}$ is independent of the basis chosen and hence is
well-defined. If $\|A\|_{2} < \infty$, then $A$ is called a
Hilbert--Schmidt operator and by this norm $B_{2}(H)$ is a Hilbert
space. Indeed, $B_2(H)$ is a special case of the Schatten
$p$-class of $H$ when $p=2$. For more details about these classes
of operators, see \cite{10,23}.

Chan \cite{9} showed that hypercyclicity can occur on the operator
algebra $B(H)$ with the strong operator topology (SOT-topology)
that is not metrizable. For example, when $T$ satisfies the
Hypercyclicity Criterion, then the left multiplication operator
$L_{T}$ is SOT-hypercyclic on $B(H)$, that is, $L_{T}$ is
hypercyclic on $B(H)$ with strong operator topology.

\section{Main results}

From now on we suppose that $H$ is a separable
infinite-dimensional Hilbert space.

\setcounter{theore}{0}
\begin{definit}$\left.\right.$\vspace{.5pc}

\noindent {\rm Let $L\hbox{:}\ B(H) \rightarrow B(H)$ be linear
and bounded. We say that $L$ is SOT-hypercyclic if there exist
some $T \in B(H)$ such that the set $Orb(L,T) = \{T, LT, L^{2}T,
\ldots\}$ is dense in $B(H)$ in the strong operator topology. Also
we say that $L\hbox{:}\ B_{2}(H) \rightarrow B_{2}(H)$ is
$\|\hbox{$\cdot$}\|_{2}$-hypercyclic if there exists some $T \in B_{2}(H)$
such that $Orb(L,T)$ is dense in $B_{2}(H)$ with
$\|\hbox{$\cdot$}\|_{2}$-topology.}
\end{definit}

\begin{definit}$\left.\right.$\vspace{.5pc}

\noindent {\rm For any operator $T \in B(H)$, define the left
multiplication operator $L_{T}\hbox{:}\ B(H) \rightarrow B(H)$ by
$L_{T}(S) = TS$ for every $S \in B(H)$.}
\end{definit}

Note that $B_{2}(H)$ is an ideal of $B(H)$ and hence
$L_{T}\hbox{:}\ B_{2}(H) \rightarrow B_{2}(H)$ is also
well-defined. We show that $B(H)$ and $B_{2}(H)$, respectively
with the strong operator topology and $\|\hbox{$\cdot$}\|_{2}$-topology,
are separable. For this, see the following Lemma~2.3.

Suppose $\{e_{i}\hbox{:}\ i \geq 1\}$ is an orthonormal basis for
a separable Hilbert space $H$ and $S(H)$ denotes the set of all
finite rank operators $T$ such that there exists $N_{T} \in
\mathbb{N}$, satisfying $Te_{i} = 0$ for $i \geq N_T$.

\begin{lem} Suppose $E = \{e_{i}\hbox{\rm :}\ i \geq 1\}$ is a basis for a
separable Hilbert space $H${\rm ,} then $S(H)$ is SOT-dense in $B(H)$
and also $\|\hbox{$\cdot$}\|_{2}$-dense in $B_{2}(H)${\rm ;} moreover{\rm
,} $S(H)$ is\break separable.
\end{lem}

\begin{proof} Suppose that $A \in B_{2}(H)$ and $\varepsilon > 0$.
Then there exist $N \in \mathbb{N}$ such that
$\sum^{\infty}_{i=N+1} \|Ae_{i}\|^{2} < \varepsilon^{2}$. Now
define the finite rank operator $F$ by $F=A$ on $[e_{k}\hbox{:}\ 1
\leq k \leq N]$ and $F=0$ on $[e_{k}\hbox{:}\  1 \leq k \leq
N]^{\bot}$. ($[e_{k}\hbox{:}\  1 \leq k \leq N]$ means the linear
span of $\{e_{k}\hbox{:}\  1 \leq k \leq N\}$). Thus
$\|A-F\|^{2}_{2} = \sum^{\infty}_{i=N+1} \|Ae_{i}\|^{2}<
\varepsilon^{2}$ and so $S(H)$ is $\|\hbox{$\cdot$}\|_{2}$-dense. Also,
of \cite{9} p.~234 implies that every $\|\hbox{$\cdot$}\|_{2}$-dense subset
of $B_{2}(H)$ is SOT-dense in $B(H)$, and so it follows that
$S(H)$ is SOT-dense. Now the proof is complete.\hfill $\Box$
\end{proof}

The following result is the main tool that we used to show that an
operator is hypercyclic. Versions of this result have appeared in
the work of Godefroy and Shapiro (\cite{13}, Theorem~1.2) and
Kitai (\cite{17}, Theorem~2.1).

\begin{propo}$\left.\right.$\vspace{.5pc}

\noindent If $T$ is a continuous operator on a separable Banach
space $X${\rm ,} then $T$ is hypercyclic if and only if for any
two non-void open sets $U$ and $V$ in $X, T^{n} U \cap V \neq
\phi$ for some positive integer $n$.
\end{propo}

Godefroy and Shapiro (\cite{13}, Corollary~1.3) also gave a
sufficient condition for hypercyclicity that is a direct
consequence of Proposition~2.4.

\begin{coro}$\left.\right.$\vspace{.5pc}

\noindent An operator $T$ on a separable Banach space $X$ is
hypercyclic if for each pair $U,V$ of non-void open subsets of
$X${\rm ,} and each neighborhood $W$ of zero in $X${\rm ,} there
are infinitely many positive integers $n$ such that both $T^{n}U
\cap W$ and $T^{n} W \cap V$ are non-empty.
\end{coro}

\begin{remar}$\left.\right.$
\begin{enumerate}
\renewcommand\labelenumi{\rm (\roman{enumi})}
\leftskip .35pc
\item In Proposition~2.4, the condition $T^{n}U \cap V \neq \phi$ is
equivalent to the condition $U \cap T^{-n}V \neq \phi$.

\item \looseness -1 If an operator $T$ is hypercyclic, then it automatically has
a dense set of hypercyclic vectors. For, if a vector $x$ is
hypercyclic for $T$, then so is $T^{n}x$ for any positive integer
$n$. Thus the condition `$T^{n} U \cap V \neq \phi$ for some
positive integer $n$', in Proposition~2.4, can be replaced by the
condition `$T^{n}U \cap V \neq \phi$ for infinitely many positive
integers $n$'.

\item Equivalent to the hypothesis of Corollary~2.5 is the
apparently weaker requirement that the sets $T^{n}U \cap W$ and
$T^{n}W \cap V$ be non-empty for a single $n$.
\end{enumerate}
\end{remar}

The following theorem shows that the converse of the above
corollary is equivalent to the Hypercyclicity Criterion. Remember
that for vectors $g,h$ in $H$ the operator $g \otimes h$ denotes a
rank one operator and is defined by $(g \otimes h)(f) = \langle f,
h \rangle g$.

\begin{theor}[\!] For any operator $T \in B(H)${\rm ,} the following
are equivalent{\rm :}
\begin{enumerate}
\renewcommand\labelenumi{\rm (\roman{enumi})}
\leftskip .15pc
\item $T$ satisfies the hypothesis of the Hypercyclicity Criterion.

\item For each pair $U,V$ of non-void open subsets of $H${\rm ,} and
each neighborhood $W$ of zero{\rm ,} $T^{n}U \cap W \neq \phi$ and
$T^{n}W \cap V \neq \phi$ for some integer $n$.
\end{enumerate}
\end{theor}

\begin{proof} It is easy to see that (i) implies (ii) (for
details see Corollary~1.4 in \cite{7}). For the converse, assume
that $T$ satisfies property (ii). First we show that for each pair
$U',V'$ of non-void $\|\hbox{$\cdot$}\|_{2}$-open subsets of $B_{2}(H)$
there is an integer $n \geq 1$ such that $U' \cap L^{-n}_{T} V'
\neq \phi$. For this, fix an orthonormal basis $E=\{e_{i}\hbox{:}\ 
i \geq 1\}$ for $H$. By using Lemma~2.3 there exist finite rank
operators $A$ and $B$ such that $A \in S(H) \cap U'$ and $B \in
S(H) \cap V'$, whence for a certain integer $N \geq 1$ we have
$A(e_{i}) = B(e_{i})=0$ for $i>N$. But for some $\varepsilon > 0$
we have
\begin{equation*}
\left\{D \in S(H)\hbox{:}\ \|D-A\|_{2} < \varepsilon\right\}
\subseteq S(H) \cap U',
\end{equation*}
and
\begin{equation*}
\left\{D \in S(H)\hbox{:}\ \|D-B\|_{2} < \varepsilon\right\}
\subseteq S(H) \cap V'.
\end{equation*}
Now consider the following open sets:
\begin{equation*}
U_{i} = \left\{ h \in H\hbox{:}\ \|h-Ae_{i}\| <
\frac{\varepsilon}{2\sqrt{N}} \right\}, V_{i} = \left\{ h \in
H\hbox{:}\ \|h-Be_{i}\| < \frac{\varepsilon}{2\sqrt{N}} \right\}
\end{equation*}

$\left.\right.$\vspace{-1.5pc}

\noindent for $i=1,2,\ldots,N$. Note that Corollary~2.5 or remark (iii)
implies that $T$ is hypercyclic. Now by using Proposition~2.4
repeatedly (indeed by remark (ii)), it follows that there exist
integers $0 = n_{0} < n_{1} \leq n_{2} \leq \cdots \leq n_{N-1}$
and $0=m_{0} < m_{1} \leq m_{2} \leq \cdots \leq m_{N-1}$ such
that
\begin{equation}
U=U_{1} \cap T^{-n_{1}} U_{2} \cap T^{-n_{2}} U_{3} \cap \cdots
\cap T^{-n_{N-1}} U_{N} \neq \phi
\end{equation}
and
\begin{equation}
V=V_{1} \cap T^{-m_{1}} V_{2} \cap T^{-m_{2}} V_{3} \cap \cdots
\cap T^{-m_{N-1}} V_{N} \neq \phi.
\end{equation}
Put $W= \{h\hbox{:}\ \|h\|< \delta\}$ where
\begin{equation}
\delta = \min\left\{\frac{\varepsilon}{2\sqrt{N} \|T\|^{n_{i-1}}},
\frac{\varepsilon}{2\sqrt{N}\|T\|^{m_{i-1}}}\hbox{:}\ i =
1,2,\ldots,N\right\}.
\end{equation}
Since $T$ satisfies the hypothesis (ii) of Theorem~2.6, then there
exists some $x \in W$ and $y \in U$ such that $T^{n}x \in V$ and
$T^{n} y \in W$ for some integer $n$. The relations (1) and (2)
imply that
\begin{equation}
\|T^{n_{i-1}} y - Ae_{i}\| < \frac{\varepsilon}{2\sqrt{N}}; \quad
\|T^{n}(T^{m_{i-1}} x) - Be_{i}\| < \frac{\varepsilon}{2\sqrt{N}}
\end{equation}
for $i=1,2,\ldots,N$. Now define $S_{1} = \sum^{N}_{i=1}
T^{n_{i-1}} y \otimes e_{i}$ and $S_{2} = \sum^{N}_{i=1}
T^{m_{i-1}} x \otimes e_{i}$. Let $S=S_{1} + S_{2}$. Then $S$ is a
Hilbert--Schmidt operator, because it has finite rank. Note that by
(3), $\|T^{m_{i-1}}x\| \leq \|T\|^{m_{i-1}} \|x\| < \delta
\|T\|^{m_{i-1}} < \frac{\varepsilon}{2\sqrt{N}}$. Now by using (4)
we get the following inequalities:
\begin{align*}
\|S-A\|_{2} &\leq \|S_{1} - A\|_{2} + \|S_{2}\|_{2}\\[.4pc]
&= \left\{ \sum^{N}_{i=1} \|S_{1} e_{i} - Ae_{i}\|^{2}\right\}^{1/2} +
\left\{ \sum^{N}_{i=1} \|S_{2} e_{i}\|^{2}\right\}^{1/2}\\[.4pc]
&= \left\{ \sum^{N}_{i=1} \|T^{n_{i-1}} y -
Ae_{i}\|^{2}\right\}^{1/2} + \left\{ \sum^{N}_{i=1}
\|T^{m_{i-1}} x\|^{2}\right\}^{1/2} < \varepsilon.
\end{align*}
Hence $S \in U'$. Also note that since $T^{n}y \in W$, by (3) we
get $\|T^{n_{i-1}}(T^{n}y)\| \leq \|T\|^{n_{i-1}} \delta <
\frac{\varepsilon}{2\sqrt{N}}$, and thus we have
\begin{align*}
\|L^{n}_{T} S-B\|_{2} &\leq \|L^{n}_{T} S_{2} - B\|_{2} +
\|L^{n}_{T} S_{1}\|_{2}\\[.4pc]
&= \left\{ \sum^{N}_{i=1} \|T^{n} S_{2} e_{i} - Be_{i}
\|^{2}\right\}^{1/2} + \left\{ \sum^{N}_{i=1} \|T^{n} S_{1}
e_{i}\|^{2} \right\}^{1/2}\\[.4pc]
&= \left\{ \sum^{N}_{i=1} \|T^{n}(T^{m_{i-1}} x) -
Be_{i}\|^{2}\right\}^{1/2}\\[.4pc]
&\quad\, + \left\{ \sum^{N}_{i=1} \|T^{n_{i-1}}
(T^{n} y)\|^{2}\right\}^{1/2} < \varepsilon.
\end{align*}
So $L^{n}_{T} S \in V'$. Now it follows that $U' \cap L^{-n}_{T}
V' \neq \phi$ and so by Proposition~2.4, $L_{T}$ is $\| \;
\|_{2}$-hypercyclic. This also implies that
$\bigoplus^{\infty}_{n=1} T\hbox{:}\ \bigoplus^{\infty}_{n=1} H
\rightarrow \bigoplus^{\infty}_{n=1} H$ is hypercyclic, because
the left multiplication operator $L_{T}\hbox{:}\ B_{2}(H)
\rightarrow B_{2}(H)$ is unitary equivalent to the operator
$\bigoplus^{\infty}_{n=1} T\hbox{:}\ \bigoplus^{\infty}_{n=1} H
\rightarrow \bigoplus^{\infty}_{n=1} H$ (see \cite{11}, p.~6). Now
Theorem~2.3 in \cite{2} implies that $T$ satisfies the
Hypercyclicity Criterion, and so the proof is now complete.\hfill
$\Box$
\end{proof}

\begin{propo}$\left.\right.$\vspace{.5pc}

\noindent If $T \in B(H)${\rm ,} then the following are
equivalent{\rm :}
\begin{enumerate}
\renewcommand\labelenumi{\rm (\roman{enumi})}
\leftskip .15pc
\item $T$ satisfies the hypothesis of the Hypercyclicity Criterion.

\item $T$ is hypercyclic and for each non-void open subset $U$ and
each neighborhood $W$ of zero{\rm ,} $T^{n}U \cap W \neq \phi$ and
$T^{-n} U \cap W \neq \phi$ for some integer $n$.
\end{enumerate}
\end{propo}

\begin{proof} By Theorem~2.6 it suffices to show that (ii) implies
(i). So let (ii) hold. By Theorem~2.6, it suffices to show that
(ii) in Theorem~2.6 holds. Since $T$ is hypercyclic, by
Proposition~2.4, $U \cap T^{-m} V \neq \phi$ for some positive
integer $m$. Let $G$ be a neighborhood of zero that is contained
in $W \cap T^{-m} W$. By condition (ii), there exists some
positive integer $n$ such that $T^{-n} G \cap (U \cap T^{-m} V)
\neq \phi$ and $G \cap T^{-n} (U \cap T^{-m} V) \neq \phi$. But
$T^{-n}G \cap (U \cap T^{-m}V)$ is a subset of $T^{-n}W \cap U$,
hence $T^{-n}W \cap U \neq \phi$. Also $G \cap T^{-n}(U \cap
T^{-m}V)$ is a subset of $T^{-m}W \cap T^{-n}(T^{-m}V) = T^{-m}(W
\cap T^{-n}V)$ which implies that $T^{-n}V \cap W \neq \phi$.
Thus, hypothesis (ii) of Theorem~2.6 holds and so the proof is
complete.\hfill $\Box$
\end{proof}

\begin{rem}{\rm We say that the sequence
$\{T_{n}\}^{\infty}_{n=}$ of bounded linear operators on a Hilbert
space $H$ is hypercyclic provided that there exists some $x \in H$
such that the collection of images $\{T_{n}x\hbox{:}\ n=1,2,\ldots\}$ is
dense in $H$. Note that Theorem~1.1, Proposition~2.4 and
Corollary~2.5 can be extended to the case where hypercyclicity of
$T$ is replaced by hypercyclicity for the sequence
$\{T_{n}\}^{\infty}_{n=1}$ of bounded linear operators that have
dense range. In particular we say that $\{T_{n}\}^{\infty}_{n=1}$
satisfies the hypothesis of the Hypercyclicity Criterion if in the
hypothesis of Theorem~1.1, we use $T_{n_{k}}$ instead of
$T^{n_{k}}$. It also implies that if the sequence
$\{T_{n}\}^{\infty}_{n=1}$ satisfies the hypothesis of the
Hypercyclicity Criterion, then $\{T_{n}\}^{\infty}_{n=1}$ is
hypercyclic (see Theorem~1.2, Corollaries~1.3 and 1.5 in
\cite{13}).}
\end{rem}

It is not difficult to see that Theorem~2.6 and Proposition~2.7
work for the sequence $\{T_{n}\}^{\infty}_{n=1}$ of bounded linear
operators provided that $T_{n}T_{m} = T_{m}T_{n}$ for each pair
$m,n$ of positive integers. Hence we can deduce the following
corollary.

\begin{coro}$\left.\right.$\vspace{.5pc}

\noindent Suppose that $\{T_{n}\}^{\infty}_{n=1}$ is a sequence of
bounded linear operators on a Hilbert space $H$ such that
$T_{n}T_{m} = T_{m}T_{n}$ for each pair $m,n$ of positive integers
and have dense range. Then the following are equivalent{\rm :}
\begin{enumerate}
\renewcommand\labelenumi{\rm (\roman{enumi})}
\leftskip .35pc
\item $\{T_{n}\}^{\infty}_{n=1}$ satisfies the hypothesis of the
Hypercyclicity Criterion.

\item For each pair $U,V$ of non-void open subsets of $H${\rm ,}
and each neighborhood $W$ of zero{\rm ,} $T_{n}U \cap W \neq \phi$
and $T_{n}W \cap V \neq \phi$ for some integer $n$.

\item $\{T_{n}\}^{\infty}_{n=1}$ is hypercyclic and for each
non-void open subset $U$ and each neighborhood $W$ of zero{\rm ,}
$T_{n}U \cap W \neq \phi$ and $T^{-1}_{n}U \cap W \neq \phi$ for
some integer $n$.\vspace{-1pc}
\end{enumerate}
\end{coro}

The following definition is introduced in \cite{2}.

\begin{definit}$\left.\right.$\vspace{.5pc}

\noindent {\rm Suppose that $T \in B(H)$ and $\{n_{k}\}$ is a
sequence of positive integers. We say that $T$ is hereditarily
hypercyclic with respect to $\{n_{k}\}$ if for any subsequence
$\{n_{k_{m}}\}$ of $\{n_{k}\}$, the sequence $\{T^{n_{k_{m}}}\}$
is hypercyclic.}
\end{definit}

Now we summarize all necessary and sufficient conditions for the
Hypercyclicity Criterion in the following corollary.

\begin{coro}$\left.\right.$\vspace{.5pc}

\noindent For any operator $T \in B(H)${\rm ,} the following are
equivalent{\rm :}

\begin{enumerate}
\renewcommand\labelenumi{\rm (\roman{enumi})}
\leftskip .35pc
\item $T$ satisfies the hypothesis of the Hypercyclicity Criterion.

\item $T$ is hereditarily hypercyclic with respect to a subsequence
$\{n_{k}\}$ of positive integers.

\item $\bigoplus^{\infty}_{i=1} T$ is hypercyclic on
$\bigoplus^{\infty}_{i=1} H$.

\item The left multiplication operator $L_{T}\hbox{\rm :}\ B_{2}(H)
\rightarrow B_{2}(H)$ is $\|\hbox{$\cdot$}\|_{2}$-hypercyclic.

\item For each pair $U,V$ of non-void open subsets of $X${\rm ,} and
each neighborhood $W$ of zero{\rm ,} $T^{n}U \cap W \neq \phi$ and
$T^{n}W \cap V \neq \phi$ for some integer $n$.
\end{enumerate}
\end{coro}

\begin{proof} The proof is an immediate consequence of Theorem~2.6,
Proposition~2.7 and Theorem~2.3 in \cite{2}.\hfill $\Box$
\end{proof}

The following proposition represents some relation between
hypercyclicity and the Hypercyclicity Criterion.

\begin{propo}$\left.\right.$\vspace{.5pc}

\noindent For any operator $T \in B(H)$ the following are
equivalent{\rm :}
\begin{enumerate}
\renewcommand\labelenumi{\rm (\roman{enumi})}
\leftskip .35pc
\item $T$ satisfies the hypothesis of the Hypercyclicity Criterion.

\item There exists a dense subset $Y$ in $X$ and a sequence
$\{n_{k}\}$ such that $\{T^{n_{k}}\}$ is hypercyclic and
$T^{n_{k}}y \rightarrow 0$ for every $y \in Y$.

\item There exists a sequence $\{n_{k}\}$ such that for each pair
$U,V$ of non-void open subsets of $H$, there is $N \geq 1$ such
that $T^{n_{k}} U \cap V \neq \phi$ for any $k \geq N$.
\end{enumerate}
\end{propo}

\begin{proof}$\left.\right.$

\noindent (i) $\rightarrow$ (ii): It follows from condition (ii)
of Corollary~2.11.

\noindent (ii) $\rightarrow$ (i): Let $T_{k} = T^{n_{k}}, U$ be
any non-void open set and also let $W$ be any open neighborhood of
zero. Then by Remark~2.8, $\{T_{k}\}_{k}$ is hypercyclic and so
there is some sequence $\{m_{k}\}$ of positive integers such that
$T_{m_{k}}W \cap U \neq \phi$ for every $k \geq 1$. Now if $y \in
U \cap Y$, then $T_{m_{k}}y = T^{n_{m_{k}}}y \rightarrow 0$ which
yields $T_{m_{k}}U \cap W \neq \phi$. It holds condition (iii) of
Corollary~2.11, hence $\{T_{k}\}$ satisfies the hypothesis of the
Hypercyclicity Criterion and so $\{T^{n_{k}}\}$ and consequently
$T$ satisfy the Hypercyclicity Criterion.

\noindent (iii) $\rightarrow$ (i): It suffices to show that
condition (iii) implies condition (v) of Corollary~2.11. For this
let $U,V$ be a pair of non-void open subsets of $H$ and $W$ be any
neighborhood of zero. Then for some integer $N$, we have
\begin{equation*}
T^{n_{k}} U \cap W \neq \phi; \quad T^{n_{k}} W \cap V \neq \phi
\end{equation*}
for any $k>N$. Thus indeed condition (v) of Corollary~2.11 is
consistent.

\noindent (i) $\rightarrow$ (iii): Note that by condition (ii)
of Corollary~2.11, $T$ is hereditarily hypercyclic with respect to
a sequence $\{n_{k}\}$ of positive integers. Now  suppose that
(iii) does not hold. So there exist some pair $U,V$ of non-void
open sets such that $T^{n_{k_{m}}} U \cap V = \phi$ for some
subsequence $\{n_{k_{m}}\}$ of $\{n_{k}\}$. But
$\{T^{n_{k_{m}}}\}$ is hypercyclic and so it is a contradiction.
Hence for every pair $U,V$ of non-void open sets, there is $N \geq
1$ such that $T^{n_{k}} U \cap V \neq \phi$ for any $k \geq N$.
The proof is now complete.\hfill$\Box$
\end{proof}

\section*{Acknowledgment}

The authors thank the referee for many interesting comments and
helpful suggestions about the paper.

\end{document}